\documentclass[12pt]{amsart}
\usepackage{amssymb}
\usepackage{amsfonts}
\usepackage{amscd}

\swapnumbers

\def\limproj{\mathop{\oalign{lim\cr
\hidewidth$\longleftarrow$\hidewidth\cr}}}

\newtheorem{theorem}[subsection]{Theorem}
\newtheorem{lem}[subsection]{Lemma}

\theoremstyle{definition}

\newtheorem{def-prop}[subsubsection]{Definition-Proposition}

\theoremstyle{remark}
\newtheorem{remark}[subsection]{Remark}

\theoremstyle{plain}

\numberwithin{equation}{subsection}

\def\boxit#1#2{\setbox1=\hbox{\kern#1{#2}\kern#1}%
\dimen1=\ht1 \advance\dimen1 by #1
\dimen2=\dp1 \advance\dimen2 by #1
\setbox1=\hbox{\vrule height\dimen1 depth\dimen2\box1\vrule}%
\setbox1=\vbox{\hrule\box1\hrule}%
\advance\dimen1 by .4pt \ht1=\dimen1
\advance\dimen2 by .4pt \dp1=\dimen2 \box1\relax}

\def\AA{{\mathbf A}}

\def\CC{{\mathbf C}}

\def\LL{{\mathbf L}}

\def\NN{{\mathbf N}}

\def\QQ{{\mathbf Q}}

\def\ZZ{{\mathbf Z}}

\def\cL{{\mathcal L}}
\def\cM{{\mathcal M}}

\def\cS{{\mathcal S}}

\def\cU{{\mathcal U}}

\def\cX{{\mathcal X}}

\def\cZ{{\mathcal Z}}
\mathchardef\alphag="7C0B
\mathchardef\betag="7C0C
\mathchardef\gammag="7C0D
\mathchardef\deltag="7C0E
\mathchardef\varepsilong="7C22
\mathchardef\varphig="7C27
\mathchardef\psig="7C20
\mathchardef\zetag="7C10
\mathchardef\epsilong="7C0F
\mathchardef\rhog="7C1A
\mathchardef\taug="7C1C
\mathchardef\upsilong="7C1D
\mathchardef\iotag="7C13
\mathchardef\thetag="7C12
\mathchardef\pig="7C19
\mathchardef\sigmag="7C1B
\mathchardef\etag="7C11
\mathchardef\omegag="7C21
\mathchardef\kappag="7C14
\mathchardef\lambdag="7C15
\mathchardef\mug="7C16
\mathchardef\xig="7C18
\mathchardef\chig="7C1F
\mathchardef\nug="7C17
\mathchardef\varthetag="7C23
\mathchardef\varpig="7C24
\mathchardef\varrhog="7C25
\mathchardef\varsigmag="7C26
\mathchardef\Omegag="7C0A
\mathchardef\Thetag="7C02
\mathchardef\Sigmag="7C06
\mathchardef\Deltag="7C01
\mathchardef\Phig="7C08
\mathchardef\Gammag="7C00
\mathchardef\Psig="7C09
\mathchardef\Lambdag="7C03
\mathchardef\Xig="7C04
\mathchardef\Pig="7C05
\mathchardef\Upsilong="7C07

\def\ord{{\rm ord}}

\begin{document}

\title[Lefschetz numbers of  monodromy]{Lefschetz numbers of iterates of the monodromy and truncated arcs}

\author{Jan Denef}
\address{University of Leuven, Department of Mathematics,
Celestijnenlaan 200B, 3001 Leu\-ven, Bel\-gium }
\email{ Jan.Denef@wis.kuleuven.ac.be}
\urladdr{http://www.wis.kuleuven.ac.be/wis/algebra/denef.html}

\author{Fran\c cois Loeser}

\address{Centre de Math{\'e}matiques,
Ecole Polytechnique,
F-91128 Palaiseau
(UMR 7640 du CNRS), {\rm and}
Institut de Math{\'e}matiques,
Universit{\'e} P. et M. Curie, Case 82,
4 place Jussieu,
F-75252 Paris Cedex 05
(UMR 7596 du CNRS)}
\email{loeser@math.polytechnique.fr}
\urladdr{http://math.polytechnique.fr/cmat/loeser/}

%\date{}
%\dedicatory{}

\begin{abstract}We express the Lefschetz number of iterates of the
monodromy of a function on a smooth complex algebraic variety 
in terms of the Euler characteristic of a space of truncated arcs. 
\end{abstract}

\maketitle

\section{Introduction}\label{int}
Let $X$ be a smooth complex algebraic variety
and let $f : X \rightarrow \CC$ be a non constant
morphism of 
complex algebraic varieties. We fix a smooth metric on $X$.
Let $x$ be a point of $f^{-1} (0)$.
We set $X_{\varepsilon, \eta}^{\times} := B (x, \varepsilon) \cap f^{-1}
(D_{\eta}^{\times})$,
with $B (x, \varepsilon) $
the open ball of radius $\varepsilon$ centered at $x$ and
$D_{\eta}^{\times} =
D_{\eta} \setminus \{0\}$,
with $D_{\eta}$ the open disk of radius $\eta$ centered at $0$.
For $0 < \eta \ll \varepsilon \ll 1$, the restriction of $f$ to
$X_{\varepsilon, \eta}^{\times}$ is a locally trivial
fibration - called the Milnor
fibration - onto 
$D_{\eta}^{\times}$ with fiber $F_{x}$, the Milnor fibre at $x$.
The action of a characteristic homeomorphism of this fibration on
cohomology
gives rise to 
the monodromy operator 
$$
M : H^{\ast} (F_{x}, \QQ)
\longrightarrow
H^{\ast} (F_{x}, \QQ).
$$
For any natural number $n$, we consider the Lefschetz number
$$
\Lambda (M^{n}) := \sum_{q \geq 0} (-1)^q \, {\rm Trace} \, [M^{n}, H^{q} (F_{x},
\QQ)],
$$
of the $n$-th iterate of $M$.
These numbers are related to 
the monodromy zeta function
$$Z (t)~: = \prod_{q \geq 0} [{\rm det} \, ({\rm Id} - tM, H^{q} (F_{x},
\QQ))]^{(-1)^{q}}$$
as follows: if one writes
$$
\Lambda (M^{n}) = \sum_{i | n} s_{i},
$$
for $n \geq 1$,
then
$$
Z (t) = \prod_{i \geq 1} (1 - t^{i})^{s_{i} / i}.
$$
In the paper \cite{A'Campo}, A'Campo gives a formula for the Lefschetz numbers
$\Lambda (M^{n})$ in terms of a resolution of the singularities of $f =
0$. The aim of the present  note is to give 
a formula for $\Lambda (M^{n})$ directly in terms of geometric objects
associated to the function $f$.
This formula will involve truncated arcs on $X$.
Let us recall from \cite{Arcs} \cite{K},
that there is a $\CC$-scheme $\cL (X)$,
the space of formal arcs on $X$, whose set of $\CC$-rational points
$\cL (X) (\CC)$ is naturally in bijection with $X (\CC [[t]])$.
Similarly, for $n \geq 0$, we can consider the space 
$\cL_{n} (X)$ of arcs modulo $t^{n +1}$: a $\CC$-rational
point of
$\cL_{n} (X) $ corresponds to a $\CC [t] / t^{n + 1} \CC [t]$-rational
point on
$X$. The space
$\cL_{n} (X)$ is canonically endowed with the structure of a complex
algebraic variety.
For instance when $X$ is the affine space $\AA^{m}_{\CC}$,
a $\CC$-rational
point of $\cL (X)$ is just an $m$-tuple of power series in the
variable $t$
with coefficients in $\CC$, while 
$\cL_{n} (X) (\CC)$ is the set of
$m$-tuples of complex polynomials of degree $\leq n$ in the
variable $t$.
Furthermore, there is a natural morphism
$$ 
\pi_{n} : \cL (X) \longrightarrow \cL_{n} (X)
$$
which corresponds to truncation on $\CC$-rational points.
Since $X$ is assumed to be smooth, the morphism $\pi_{n}$ is
surjective
(also for $\CC$-rational points).
In what follows, we identify $\cL (X)$ with its set of 
$\CC$-rational points, and similarly  for
$\cL_n (X)$. With a ``complex algebraic variety'', we mean the set of
$\CC$-rational points of a separated reduced scheme of finite type
over $\CC$, not necessarly irreducible.
Now, we consider the set $\cX_{n}$ 
of points
$\varphi$ in $\cL_{n} (X)$ such that $\pi_{0} (\varphi) = x$ and
such that the $t$-valuation of $f (\varphi)$ is exactly $n$.
This subset $\cX_{n}$ is 
a locally closed subvariety of
$\cL_{n} (X)$ and we may consider the morphism
$$
\bar f_{n} : \cX_{n} \longrightarrow \CC^{\times},
$$
sending a point $\varphi$ in $\cX_{n}$ to the coefficient of
$t^{n}$ 
in $f (\varphi)$.
There is a natural action of $\CC^{\times}$ on $\cX_{n}$
given by $a \cdot \varphi (t) = \varphi (a t)$. Since
$\bar f_{n} (a \cdot \varphi) =
a^{n} \bar f_{n} (\varphi)$, it follows
that $\bar f_{n}$ is a locally trivial fibration and has a geometric
monodromy of order $n$.
We denote by  $\cX_{n, 1}$ the fiber $\bar f_{n}^{-1} (1)$.
The group $\mu_n$ of $n$-th roots of unity in $\CC$
acts naturally on $\cX_{n, 1}$, by restriction of the
$\CC^{\times}$-action on $\cX_n$. 

\medskip

We can now state the main result of this note.

\begin{theorem}\label{MT}For every integer $n \geq 1$, the Lefschetz number 
$\Lambda (M^{n})$ is equal to $\chi (\cX_{n, 1})$, the
Euler
characteristic of $\cX_{n, 1} $.
\end{theorem}

In fact, we shall deduce Theorem \ref{MT} from a more general result,
Theorem \ref{PT}, where we give a formula for the class of
$\cX_{n, 1}$ in the ring $\cM_{\rm loc}$ obtained
by localisation of the class of the affine line in the Grothendieck ring of
complex algebraic varieties, whose definition is recalled at the
beginning of section \ref{next}.
Actually, it would
be also possible to prove Theorem \ref{MT}
as an easy consequence of Theorem 2.2.1 in \cite{Motivic}, by a reasonning
involving motives and Hodge polynomials.
In fact, Theorem \ref{PT} may be further generalized
to take in account the monodromy action. This is done in 
Theorem \ref{MTmono}, once we introduced the ``monodromic'' Grothendieck ring
$\cM^{\rm mon}_{\rm loc}$.

\medskip

We denote by $T_n$ the monodromy operator
$$
T_n : H^{\ast} (\cX_{n,1}, \QQ)\medskip
\noindent{\bf }
\medskip
\longrightarrow
H^{\ast} (\cX_{n,1}, \QQ)
$$
of the locally trivial fibration $\bar f_{n}~: \cX_n \rightarrow \CC^{\times}$.
Note that $T_n$ is induced by the automorphism
$\varphi (t) \mapsto \varphi (e^{2 \pi i /n} t)$ of
$\cX_{n,1}$.
For any $d$ in
$\NN$, we denote the Lefschetz number of
$T_n^d$ by 
$$
\Lambda (T^{d}_n) := \sum_{q \geq 0} (-1)^q \, {\rm Trace} \, [T^{d}_n, H^{q} (\cX_{n,1},
\QQ)].
$$
Since $T^n_n$
is the identity, we have
$\Lambda (T^{d}_n) = \Lambda (T^{\gcd(d,n)}_n)$.
Hence to know the
$\Lambda (T^{d}_n)$'s for every $d$ in $\NN$, it is enough to know them
for every $d$
dividing $n$. This information is provided by the following result:

\begin{theorem}\label{sec}If $n \geq 1$ and $d$ divides $n$, then
$\Lambda (T^{d}_n) = \Lambda (M^d)$.
\end{theorem}

Finally, in section \ref{further}, we extend Theorem \ref{MTmono}
to the case of quasi-projective varieties over a field of characteristic zero.
This enables us to construct  what we believe to
be the ``virtual motivic incarnation'' of the
the Milnor fibre at $x$ in
$\cM^{\rm mon}_{k, \rm loc}$, the analogue over $k$ of the ring
$\cM^{\rm mon}_{\rm loc}$. Taking Euler characteristic with values
into
virtual Chow motives,
we apply this to settle an issue that remained open in \cite{Motivic}.

\medskip

The first author is endebted to Eduard Looijenga for a conversation
in November 1998 which influenced the present paper.
The second author would like to thank Paul Seidel for interesting
discussions.

\section{Calculation of a motivic volume and proof of
the main results}\label{next}
To state Theorem \ref{PT}, we shall start with some reminders from motivic
integration as developped in \cite{Arcs} and \cite{MK}.

\subsection{}We denote by $\cM$ the abelian group
generated by symbols $[S]$, for $S$ a complex
algebraic variety, with the relations $[S] = [S']$ if $S$ and $S'$ are
isomorphic and $[S] = [S'] + [S \setminus S']$
if $S'$ is Zariski
closed in $S$. There is a natural ring structure
on $\cM$, the product being
induced by the cartesian product of varieties, and
to any constructible set $S$ in some complex algebraic
variety  one naturally associates
a class $[S]$ in $\cM$.
We denote by $\LL$ the class
of the affine line in $\cM$ and
we denote
by $\cM_{\rm loc}$ the localisation 
$\cM_{\rm loc} :=  \cM [\LL^{-1}]$.

\medskip
Let $X$ be a smooth complex algebraic variety of pure dimension $m$.
We call a subset $A$ of $\cL (X)$ cylindrical at level $n$
if $A = \pi_{n}^{-1} (C)$, with $C$ a constructible subset of
$\cL_{n} (X)$. We say that $A$ is cylindrical if it is cylindrical
at some level $n$.

Let $X$, $Y$ and $F$ be complex algebraic varieties,
and let
$A$, resp. $B$, be a constructible subset of $X$,
resp. $Y$. We say that
a map
$\pi : A \rightarrow B$ is a piecewise morphism if
there exists a finite partition of
the domain of $\pi$ into locally closed subvarieties of $X$
such that the restriction of $\pi$ to any of these
subvarieties is a morphism of
varieties.
We say that
a map
$\pi : A \rightarrow B$ is a
piecewise trivial fibration with fiber
$F$, if there exists a finite partition of $B$ in subsets $S$ which are
locally closed subvarieties of
$Y$ such that $\pi^{- 1} (S)$ is a
locally closed subvariety of  $X$ and
isomorphic, as a complex algebraic
variety, to $S \times F$, with $\pi$
corresponding under the isomorphism to the projection
$S \times F \rightarrow S$. We say that the map $\pi$ is
a
piecewise trivial fibration over some constructible subset $C$ of
$B$,
if the restriction of $\pi$ to $\pi^{- 1} (C)$ is a piecewise 
trivial fibration onto $C$.

We call a cylindrical subset $A$ of $\cL (X)$
stable at level $n \in \NN$ if $A$ is cylindrical
at level $n$ and $\pi_{k + 1} (\cL (X)) \rightarrow \pi_{k} (\cL (X))$
is a piecewise
trivial fibration over $\pi_{k} (A)$ with fiber $\AA^{m}_{\CC}$, for all $k
\geq n$.
We call $A$ stable if it is stable at some level $n$.

Since $X$ is smooth, any cylindrical subset $A$ of $\cL (X)$ is stable
(at the same level),
by Lemma 4.1 of \cite{Arcs}.
Denote by $\CC_{0}$ the family of stable cylindrical subsets of $\cL
(X)$.
Clearly there exists a unique additive measure
$$
\tilde \mu : \CC_{0} \longrightarrow \cM_{\rm loc}
$$
satisfying
$$
\tilde \mu (A) = [\pi_{n} (A)] \, \LL^{- (n + 1) m}
$$
when $A \in \CC_{0}$ is stable at level $n$.

In particular, the relation
$$
[\cX_{n, 1}] = \tilde \mu (\cZ_{n, 1}) \, \LL^{(n + 1) m}
$$
holds in $\cM_{\rm loc}$, where $\cZ_{n, 1}$ is the set
of points
$\varphi$ in $\cL (X)$ such that $\pi_{0} (\varphi) = x$,
such that the $t$-valuation of $f (\varphi)$ is exactly $n$, and such
that the coefficient of $t^{n}$ in $f (\varphi)$ is equal to 1.

The following geometric lemma, which is a special case of Lemma 3.4 in
\cite{Arcs}, will play a crucial role
in the proof of Theorem \ref{PT}.

\begin{lem}\label{fib}Let $X$ and $Y$ be connected smooth complex
algebraic
varieties
and let $h : Y \rightarrow X$ be a birational morphism.
For $e$ in $\NN$, let $\Delta_e$ be the subset
of $\cL (Y)$ defined
by
$$
\Delta_e 
:= \{\varphi \in Y (\CC[[t]]) \bigm \vert {\rm ord}_t {\rm det} \,
{\rm Jac}_{h} (\varphi)
= e \},
$$
where ${\rm Jac}_{h} (\varphi)$ is the jacobian of $h$ at $\varphi$.
For $k$ in $\NN$, let
$h_{k \ast} : \cL_k (Y) \rightarrow \cL_k (X)$ be the morphism
induced by $h$, and let $\Delta_{e, k}$ be the image of
$\Delta_e$ in $\cL_k (Y)$. If $k \geq 2e$, the following holds.
\begin{enumerate}
\item[a)]The constructible subset $\Delta_{e, k}$ of
$\cL_k (Y)$ is  a union of fibers of
$h_{k \ast}$.
\item[b)]The restriction of $h_{k \ast}$ to $\Delta_{e, k}$
is a piecewise trivial fibration
with fiber
$\AA^e_{\CC}$ onto its image.
\end{enumerate}
\end{lem}

\subsection{}\label{resol}Now we shall use Lemma \ref{fib}
to compute $\tilde \mu (\cZ_{n, 1})$
on a resolution of $f$.
Let $D$ be the divisor defined by $f = 0$ in
$X$. Let $(Y, h)$ be a resolution of $f$. By this, we mean that 
$Y$ is a smooth and connected complex
algebraic variety, $h : Y \rightarrow
X$
is proper, that the restriction
$h : Y \setminus h^{-1} (D) \rightarrow
X \setminus D$ is an isomorphism, and that
$(h^{-1} (D))_{\rm red}$ has only normal crossings as a
subvariety
of $Y$. Furthermore,
we choose $h$ in
such a way that 
$(h^{-1} (x))_{\rm red}$ is a union of
irreducible (smooth) components
of
$(h^{-1} (D))_{\rm red}$, which we shall denote by
$E_{i}$, $i \in J$. For each $i \in J$, denote by
$N_{i}$ the multiplicity of $E_{i}$ in the divisor of
$f \circ h$ on $Y$, and by
$\nu_{i} - 1$ the multiplicity of $E_{i}$ in the divisor
of $h^{\ast} dx$, where $d x$ is a local non vanishing volume form at $x$,
{\it i.e.} a local generator of the sheaf of differential forms of
maximal degree at $x$.
For $i \in J$ and $I \subset J$, we consider the varieties
$E_{i}^{\circ} := E_{i} \setminus \cup_{j \not= i} E_{j}$,
$E_{I} := \cap_{i \in I} E_{i}$, and
$E_{I}^{\circ} := E_{I} \setminus \cup_{j \in J \setminus I} E_{j}$.
We shall also set $m_{I} = \gcd (N_{i})_{i \in I}$.
We introduce an unramified Galois cover $\widetilde E_{I}^{\circ}$ of
$E_{I}^{\circ}$, with Galois group $\mu_{m_I}$, as follows.
Let $U$ be an affine Zariski open subset of $Y$, such that, on $U$,
$f \circ h = u v^{m_I}$, with $u$ a unit on $U$ and $v$ a morphism
form $U$ to $\AA^1_{\CC}$. Then the restriction of 
$\widetilde E_{I}^{\circ}$ above $E_{I}^{\circ} \cap U$,
denoted by 
$\widetilde E_{I}^{\circ} \cap U$,
is defined as
$$\Bigl\{(z, y) \in \AA^1_{\CC} \times
(E_{I}^{\circ} \cap U)
\Bigm | z^{m_I} = u^{-1}
\Bigr\}.
$$
Note that $E_{I}^{\circ}$ can be covered by such affine open subsets
$U$ of $Y$.
Gluing together the covers 
$\widetilde E_{I}^{\circ} \cap U$, in the obvious way (cf. the proof of
Lemma 3.2.2 in \cite{Motivic}), we obtain the cover
$\widetilde E_{I}^{\circ}$ of
$E_{I}^{\circ}$ which has a natural
$\mu_{m_I}$-action (obtained by multiplying the $z$-coordinate
with the elements of $\mu_{m_I}$).

\begin{theorem}\label{PT}With the previous notations, the following
relation holds in
$\cM_{\rm loc}$:
\begin{equation}\label{ME}
[\cX_{n, 1}]
= \LL^{nm} \,
\sum_{I \subset J \atop I \not= \emptyset} 
(\LL - 1)^{|I| - 1} [\widetilde E_{I}^{\circ}]
\Biggl ( \sum_{
k_{i} \geq 1, i \in I \atop
\sum k_{i} N_{i} = n}
\LL^{ - \sum_{i \in I} k_{i} \nu_{i}}\Biggr ).
\end{equation}
\end{theorem}

\begin{proof}Let
$\widetilde \cZ_{n, 1}$ be the preimage of $\cZ_{n, 1}$ in
$\cL (Y)$. Remark that, the morphism
$h$ being proper, the induced function
$h_{\ast} : \widetilde \cZ_{n, 1} \rightarrow \cZ_{n, 1}$
is bijective. Now, for $e \geq 0$, define
$\widetilde \cZ_{n, 1, e} $ as the set of points $\varphi$
in
$\widetilde \cZ_{n, 1}$ such that 
${\rm ord}_t {\rm det} \, {\rm Jac}_{h} (\varphi)
= e$. 

\begin{lem}\label{triv}Let $I$ be a non empty subset of $J$.
Let $U$ be an 
affine Zariski
open subset of $Y$, such that, on $U$,
$f \circ h = u \prod_{i \in I} y_{i}^{N_{i}}$ and
${\rm det} \, {\rm Jac}_{h} = v \prod_{i \in I} y_{i}^{\nu_{i} - 1}$,
with $u$ and $v$ units on $U$ and $y_{i}$ a regular
function on $U$ with divisor
$E_{i} \cap U$. Let $k_{i}$, $i \in I$, be natural numbers with $\sum_{i
\in I} k_{i} N_{i} = n$, $k_{i} \geq 1$.
Let $\cU_{(k_{i})}$ be the set 
of points $\varphi$
in
$\cL_{n} (U)$ such that $\ord_{t} y_{i} (\varphi) = k_{i}$, for $i
\in I$, and $\bar f_{n} (h_{n \ast} (\varphi)) = 1$,
where $h_{n \ast} :  \cL_n (Y) \rightarrow \cL_n (X)$
is the morphism induced by $h$.
Then
$$
[\cU_{(k_{i})}] = (\LL - 1)^{|I| - 1} [\widetilde E_{I}^{\circ} \cap U]
\, \LL^{mn - \sum_{i \in I}
k_{i}}
$$
in $\cM_{\rm loc}$.
\end{lem}

\begin{proof}Note that
the projection $\cL_n
(U) \rightarrow U$ maps $\cU_{(k_i)}$ into $E_{I}^{\circ} \cap U$.
Thus we may assume that $E_{I}^{\circ} \cap U$ is not empty, and - using
additivity - that there exist $m - |I|$ functions on $U$
which, together with the functions $(y_i)_{i \in I}$,
induce an {\'e}tale map $U \rightarrow \AA_{\CC}^m$.
By Lemma 4.2 of \cite{Arcs} (with $n = e = 0$), this map induces an
isomorphism
$\cL_n (U) \simeq U \times_{\AA_{\CC}^m} \cL_n (\AA_{\CC}^m)$.
It follows now from the very definitions that
$$[\cU_{(k_{i})}] = [W_{I}]  \, \LL^{mn - \sum_{i \in I}k_{i}},$$ with
$$W_{I} := 
\Bigl \{(z_{i}, y) \in  (\CC^{\times})^{|I|} \times
(E_{I}^{\circ} \cap U)
\Bigm |
\prod_{i \in I} z_{i}^{N_{i}} u = 1
\Bigr \},$$ and we have $$[W_{I}] = 
(\LL - 1)^{|I| - 1}  \, [\widetilde E_{I}^{\circ} \cap U].$$
To verify the last equality, consider an automorphism
$z \mapsto (z^{a_i})_{i \in I}$ of $(\CC^{\times})^{|I|}$,
with $(a_i)_{i \in I}$ a basis for $\ZZ^{|I|}$ and $a_1 = (N_i /
m_I)_{i \in I}$.
\end{proof}

By covering $Y$ with affine Zariski
open subsets $U$  verifying the assumptions in Lemma \ref{triv},
one sees
that  all the subsets 
$\widetilde \cZ_{n, 1, e} $ of $\cL (Y)$ are cylindrical
and that there exists $e_{0}$ such that
$\widetilde \cZ_{n, 1, e} $ is empty for
$e > e_{0}$.
Now set 
$\cZ_{n, 1, e} = h_{\ast} (\widetilde \cZ_{n, 1, e})$. Since
$\pi_{k} \circ h_{\ast} = h_{k \ast} \circ \pi_{k}$, with the notation
of Lemma \ref{fib}, it follows from assertion a) of Lemma \ref{fib}
that
the subsets $\cZ_{n, 1, e} $ of $\cL (X)$ are cylindrical.
Since $\cZ_{n, 1} $ is equal to the disjoint union of the subsets
$\cZ_{n, 1, e} $ for $e \leq e_{0}$,
we have
\begin{equation*}
\tilde \mu (\cZ_{n, 1})
=
\sum_{e \leq e_{0}}\tilde \mu (\cZ_{n, 1, e}).
\end{equation*}
Now it follows from Lemma \ref{fib}
that
\begin{equation*}
\tilde \mu (\cZ_{n, 1, e})
=
\LL^{-e} \, \tilde \mu (\widetilde \cZ_{n, 1, e}).
\end{equation*}
By using Lemma \ref{triv}, one gets, for every $U$ as in \ref{triv},
\begin{equation*}
\tilde \mu (\widetilde \cZ_{n, 1, e} \cap \cL (U))
= \LL^{-m}
\sum_{I \subset J \atop I \not= \emptyset}  (\LL - 1)^{|I| - 1}
\,
[\widetilde E_{I}^{\circ} \cap U]
\sum_{{k_{i} \geq 1, i \in I \atop
\sum k_{i} N_{i} = n, \sum k_{i} (\nu_{i} - 1) = e}}
\LL^{ - \sum_{i \in I}  k_{i}},
\end{equation*}
and the result follows by additivity of $\tilde \mu$.
\end{proof}

\subsection{Proof of Theorem \ref{MT}}We use a resolution 
$(Y, h)$ of $f$ satisfying the conditions in \ref{resol}.
We shall view the Euler
characteristic
of complex constructible sets as a ring morphism
$\chi : \cM \rightarrow \ZZ$.
Since $\chi (\LL) = 1$, this morphism
extends uniquely to a ring morphism
$\chi : \cM_{\rm loc} \rightarrow \ZZ$.
Since $\chi ((\LL - 1)^{|I| - 1}) = 0$
when $|I| > 1 $, it follows from Theorem \ref{PT}
that
$$
\chi (\cX_{n, 1}) = \sum_{N_{i} \vert n}
N_{i} \,\chi (E_{i}^{\circ}).
$$
The result follows since $\Lambda (M^{n})$ is equal to the right hand
side of the previous formula by A'Campo's formula for $\Lambda (M^{n})$
given in
\cite{A'Campo}. \hfill \qed

\begin{remark}As observed by Paul Seidel, Theorem \ref{MT}
bears some similarity with properties of
Floer homology for a symplectic automorphism
(see, \textit{e.g.},
\cite{DS}).
\end{remark}

\begin{remark}\label{without}Of course, it is also possible to
prove Theorem \ref{MT} directly,
without considering the ring $\cM_{\rm loc}$,
by using Lemma \ref{fib} and working with Euler characteristics
all the way.
\end{remark}

\subsection{}\label{mono}
To take in account the monodromy action we introduce a
ring
$\cM^{\rm mon}$.
As an abelian group $\cM^{\rm mon}$ is generated by symbols
$[S, \tau]$, with $S$ a complex
algebraic variety and $\tau$
an automorphism of $S$.
The relations are
$[S, \tau] = [S', \tau']$ if the pairs
$(S, \tau)$ and $(S', \tau')$ are
isomorphic, $[S, \tau] = [S', \tau_{\vert S'}] + [S \setminus S',
\tau_{\vert S \setminus S'}]$
when $S'$ is Zariski
closed in $S$ and stable under $\tau$,
and
$[S \times \AA^n_{\CC}, \sigma] = [S \times \AA^n_{\CC}, \sigma']$
whenever
$\sigma$ and $\sigma'$ are liftings\footnote{meaning that
$\tau \circ p = 
p \circ \sigma = p \circ \sigma'$, where $p$ is the projection
of $S \times \AA^n_{\CC}$ onto $S$.} of a same automorphism $\tau$
of $S$.
There is a natural ring structure
on $\cM^{\rm mon}$, 
induced by the cartesian product of varieties.
We denote by $\LL$ the class of $(\AA^1_{\CC}, {\rm id})$ in
$\cM^{\rm mon}$
and
we denote
by $\cM^{\rm mon}_{\rm loc}$ the localisation 
$\cM^{\rm mon}_{\rm loc} :=  \cM^{\rm mon} [\LL^{-1}]$.

We will often write $[S]$ instead of $[S, \tau]$ when it is clear from
the the context what $\tau$ is. For example we write
$[\cX_{n, 1}]$ and $[\widetilde E_{I}^{\circ}]$, the automorphism
being the obvious one induced by multiplication by
$e^{2 \pi i /n}$, resp. $e^{2 \pi i /m_I}$.

Let $T$ be an automorphism of the scheme
$\cL (X)$ which permutes the fibers of $\pi_n$ for every $n$.
Denote by $\CC_{0, T}$ the family of stable cylindrical subsets $A$
of $\cL
(X)$ with $T (A) = A$.
Clearly there exists a unique additive measure
$$
\tilde \mu^{\rm mon} : \CC_{0, T} \longrightarrow \cM^{\rm mon}_{\rm loc}
$$
satisfying
$$
\tilde \mu^{\rm mon} (A) = [\pi_{n} (A)] \, \LL^{- (n + 1) m}
$$
when $A \in \CC_{0, T}$ is stable at level $n$.

\begin{theorem}\label{MTmono}Relation (\ref{ME})
holds in $\cM^{\rm mon}_{\rm loc}$.
\end{theorem}

\begin{proof}The proof of Theorem \ref{PT} carries over literally to
the monodromic situation.
\end{proof}

\subsection{}\label{2.8}For any complex algebraic variety with an action of a
finite abelian group
$G$, and for any character $\alpha$ of $G$,
we denote by $H^{\ast} (X, \CC)_{\alpha}$ the part of $H^{\ast} (X,
\CC)$
on which $G$ acts by multiplication by $\alpha$, and we set
$$
\chi (X, \alpha) :=
\sum_{q \geq 0} (-1)^q {\rm dim} \, H^q (X, \CC)_{\alpha}.
$$

\subsection{Proof of Theorem \ref{sec}}The monodromy zeta
function $Z (t)$ of the Milnor fibration being equal to
$\prod_{i \in J} (1 - t^{N_i})^{\chi (E_{i}^{\circ})}$
by \cite{A'Campo}, Theorem \ref{sec} is equivalent to the assertion
that the monodromy zeta
function of the fibration $\bar f_n$ is equal to
$$\prod_{i \in J \atop N_i | n} (1 - t^{N_i})^{\chi (E_{i}^{\circ})}.$$
But this assertion is equivalent to the validity of the equality
\begin{equation*}
\chi (\cX_{n, 1}, \alpha) =
\sum_{i \in J \atop {\rm ord} (\alpha) | N_i | n}
\chi (E_{i}^{\circ}),
\end{equation*}
for every character $\alpha$ of $\mu_n$, which is
a direct consequence of Theorem
\ref{MTmono}. \hfill \qed

\begin{remark}In fact, Theorem \ref{MTmono}
still remains valid if one changes
the definition of
$\cM^{\rm mon}_{\rm loc}$ in \ref{mono}
by imposing the relation
$[S \times \AA^n_{\CC}, \sigma] = [S \times \AA^n_{\CC}, \sigma']$
only when
$\sigma$ and $\sigma'$
are cartesian products which coincide on $S$. To verify this claim
one needs a straighforward refinement of Lemma 4.1 of \cite{Arcs} and of
Lemma \ref{fib}.
\end{remark}

\section{Some further results}\label{further}
\subsection{}There is an algebraic analogue of Theorem \ref{MTmono}
over an arbitrary field $k$ of characteristic zero.
To state it, we first define the ring $\cM_{k, \rm loc}^{\rm mon}$,
which is the algebraic analogue of the ring $\cM^{\rm mon}_{\rm loc}$.
For $n \geq 1$ an integer, we denote by $\mu_n$ the group {\emph
{scheme}}
over $k$ of $n$-th roots of unity.
Note that it is not assumed that all geometric points of 
$\mu_n$ are rational over $k$.
By an action of $\mu_n$ on a quasi-projective scheme over $k$,
we mean an action in the sense of group schemes and schemes over $k$.
Set $\hat \mu := \limproj_n \mu_n$.
By an action of $\hat \mu$ on a quasi-projective scheme over $k$,
we mean an action which factors through a suitable 
$\mu_n$. 
We define $\cM^{\rm mon}_{k, \rm loc}$
in the same way as $\cM^{\rm mon}_{\rm loc}$, working now with pairs consisting
of a quasi-projective scheme over $k$ and a $\hat \mu$-action on it.

\medskip

Let $X$ be a smooth algebraic variety over $k$ and let $f~: X
\rightarrow
\AA^1_{k}$ be a non constant morphism.
Similarly as in the complex case one defines the arc space $\cL (X)$
and the quasi-projective variety $\cX_{n, 1}$ with a natural $\hat
\mu$-action.
One defines also similarly as before
resolutions $(Y, h)$ of $f$ and the varieties
$\widetilde E_I^{\circ}$ with
$\hat
\mu$-action.

The proof of Theorem \ref{MTmono} carries over
to the algebraic case to give the following:

\begin{theorem}\label{MTmonoalg}Let $X$ be a smooth algebraic variety over $k$ and let $f~: X
\rightarrow
\AA^1_{k}$ be a non constant morphism.
Let $(Y, h)$ be a resolution of $f$.
With the previous notations, 
relation (\ref{ME})
holds in
$\cM^{\rm mon}_{k, \rm loc}$.
\end{theorem}

\subsection{}Consider the power series
$$
P (T)
:= 
\sum_{ n \geq 1}
[\cX_{n, 1}]
\,\LL^{- nm}\,
T^n \,
$$
in the variable $T$ over the ring $\cM^{\rm mon}_{k, \rm loc}$.
It follows directly from Theorem \ref{MTmonoalg}
that $P (T)$ is a rational power series and that its
``limit for $T \rightarrow \infty$'' (cf. section 4 of
\cite{Motivic}) is equal to $-\cS$, where
\begin{equation}\label{mil}
\cS :=
\sum_{\emptyset \not= I \subset J}
(1 - \LL)^{|I| - 1} \,
[\widetilde E_{I}^{\circ}].
\end{equation}
In particular it follows that the right hand side of
(\ref{mil}) is independent
of the resolution $(Y,h)$, as an element of $\cM^{\rm mon}_{k, \rm loc}$.
As we shall explain in \ref{last}, we believe that
$\cS$ is the  the ``virtual motivic incarnation'' of 
the Milnor fibre at $x$.

\subsection{}We denote by $E$ the smallest
subfield of $\CC$ containing all 
roots of unity. Assume $E$ is contained in $k$.
To any quasi-projective variety $X$ over $k$
with a $\hat \mu$-action, and to any character $\alpha$ of $\hat \mu$
of finite order,
we associate, as in Theorem 1.3.1 of \cite{Motivic} (see also
\cite{TS}), an element 
$\chi_{{\rm mot}, c} (X, \alpha)$ of the Grothendieck group
$K_0 ({\rm Mot}_{k, E})$ of the category
${\rm Mot}_{k, E}$ of Chow motives over $k$ with coefficients in
$E$. One can check that $\chi_{{\rm mot}, c} (\_, \alpha)$ 
factorizes through a ring morphism
$\cM^{\rm mon}_{k, \rm loc} \rightarrow K_0 ({\rm Mot}_{k, E})$.
One verifies that $\chi_{{\rm mot}, c} (\_, \alpha)$  respects
the last relation in the definition of 
$\cM^{\rm mon}_{k, \rm loc}$ by an argument similar to 
the proof of Proposition 2.6 in
\cite{MK}.

We  still write $\LL$ for $\chi_{{\rm mot}, c} (\AA_{k}^1, 1)$.
When $k = \CC$, the topological Euler characteristic of
$\chi_{{\rm mot}, c} (X, \alpha)$ is equal to
$\chi (X, \alpha)$, with the notation of \ref{2.8}.

\subsection{}\label{last}Let $\alpha$ be a character of $\widehat \mu$
of
order $d$ and
set $S_{\alpha} := \chi_{{\rm mot}, c} (\cS, \alpha^{-1})$. It follows from
the definition that
$$
S_{\alpha} =
\sum_{\emptyset \not= I \subset J \atop d | m_I}
(1 - \LL)^{|I| - 1} \,
\chi_{{\rm mot}, c} (\widetilde E_{I}^{\circ}, \alpha^{-1}).
$$
This element $S_{\alpha}$ plays a key role in section 4 of
\cite{Motivic},
where it was proved that modulo 
$(\LL - 1)$-torsion $S_{\alpha}$ does not depend on the chosen
resolution $(Y, h)$ of $f$.
However it follows now from the above considerations that
$S_{\alpha}$, as an element of 
$K_0 ({\rm Mot}_{k, E})$, does not depend on the chosen resolution.
As explained in \cite{Motivic}, we believe that $S_{\alpha}$ is
the ``virtual motivic incarnation'' of the $\alpha$-isotypic
part of the Milnor fibre. It was shown in \cite{Motivic} that
this is
indeed true at the level of  $\CC$-Hodge realizations.

\bibliographystyle{amsplain}

\end{document}